\documentstyle[12pt]{article}

\newcommand{\braid}[3]{{#1}$\lower4pt\hbox{$\oo\atop\raise4pt
           \hbox{$\scriptscriptstyle {#3} $}$}${#2}}
\newcommand{\twist}[3]{{#1}${\,\scriptscriptstyle {#3}}\atop\raise9pt
           \hbox{$\scriptstyle\oo$} ${#2}}

\newcommand{\be}{\begin{eqnarray}}
\newcommand{\ee}{\end{eqnarray}}
\newcommand{\n}{\nonumber }
\newcommand{\oo}{\otimes  }

\newcommand{\si}{\sigma}

\newcommand{\x}{\xi}

\newcommand{\tX}{\tilde X}
\newcommand{\al}{\alpha}

\begin{document}

\begin{titlepage}
\begin{center}
{\Large\bf  Extended jordanian twists for Lie algebras}
\end{center}
\vskip 1cm
\centerline{\large \bf
P.P. Kulish\footnote{Supported by the RFFI grants NN 96-01-00851
and 98-01-00310.}}
\vskip 1mm
\begin{center}
{St.Petersburg
Department of the Steklov
Mathematical Institute,}\\
{Fontanka 27, St.Petersburg, 191011,
Russia }
\end{center}
\vskip 5mm
\centerline{\large \bf
V.D. Lyakhovsky\footnote{Supported by the RFFI grant N 97-01-01152}}
\vskip 1mm
\begin{center}
{Department of Physics}\\
{St.Petersburg
State University}\\
{Ulianovskaya 1, Petrodvorets, St.Petersburg, 198904,
Russia }
\end{center}
\vskip 5mm
\centerline{\large \bf
A.I. Mudrov}
\vskip 1mm
\begin{center}
{Institute of Physics}\\
{St.Petersburg State University}\\
{Ulianovskaya 1, Petrodvorets, St.Petersburg, 198904,
Russia }
\end{center}
\vskip 5mm

\begin{abstract}
Jordanian quantizations of Lie algebras are studied using the factorizable
twists. For a restricted Borel subalgebras ${\bf B}^{\vee}$ of $sl(N)$
the explicit expressions are obtained for the twist element ${\cal F}$,
universal ${\cal R}$-matrix and the corresponding canonical element
${\cal T}$. It is shown that the twisted Hopf algebra ${\cal U}_{\cal F}
({\bf B}^{\vee})$ is self dual. The cohomological properties of the
involved Lie bialgebras are studied to justify the existence of a
contraction from the Dinfeld-Jimbo quantization to the jordanian one.
The construction of the twist is generalized to a certain type of
inhomogenious Lie algebras.

\end{abstract}

\end{titlepage}

\section{Introduction}

The thorough formulation of the theory of quantum groups by Drinfeld \cite
{DRIN} includes two types of Hopf algebras: triangular (with the universal $%
{\cal R}$-matrix satisfying the relation ${\cal R}_{21} {\cal R}=1$) and
quasitriangular (with ${\cal R}_{21} {\cal R}\not=1$). Deformations of
universal enveloping of simple Lie algebras initiated by the quantum inverse
scattering method and discovered by Drinfeld and Jimbo \cite{DRIN,JIMB}
belong to the latter class. In the framework of the deformation quantization
theory \cite{FLA} these quantum algebras correspond to Lie bialgebras with
classical $r$-matrix
$$
r_{DJ}=\sum_{i=1}^{k}t_{ij}H_i\otimes H_j + \sum_{\alpha\in
\Phi_+}E_\al\otimes E_{-\alpha},
$$
where $k$ is the rank, $t_{ij}$ is the inverse Cartan matrix, and $\Phi_+$
is the set of positive roots. This $r_{DJ}$ is one of the multitude of
solutions to the classical Yang-Baxter equation. The detailed classification
of solutions was performed for simple Lie algebras in \cite{BEL}. Only for
some of these classical $r$-matrices the corresponding quantum $R$-matrices
are known explicitly.

Although the existence of quantization for any Lie bialgebra is now proved
\cite{ETI}, the explicit knowledge of $R$-matrix as an algebraic element $%
{\cal R}$ or a matrix in some irreducible representations is required in the
FRT approach \cite{FAD} and in variety of applications of quantum groups.
One can mention the universal $R$-matrix of the quantum algebra ${\cal U}%
_q(sl(2))$ \cite{DRIN} which is a building block for the universal $R$%
-matrices for other simple Lie and Kac-Moody algebras. As about triangular
quantum groups and twisting \cite{D2,D3}, 
the well known example is the jordanian quantization of 
$sl(2)$ or, more exactly, of its Borel
subalgebra ${\bf B}_+$ ($\{h,x|[h,x]=2x\}$)
with $r=h\otimes x-x\otimes h=h\wedge x$ \cite{DRIN} and the triangular $R$%
-matrix ${\cal R}={\cal F}_{21}{\cal F}^{-1}$ defined by the twisting
element \cite{GER,OGIEV}
\begin{equation}
\label{og-twist}{\cal F}=\exp \{\frac 12h\otimes \ln (1+2\xi x)\}.
\end{equation}
This quantum algebra ${\cal U}_\xi (sl(2))$ also found numerous applications
from the deformed Heisenberg $XXX$-spin chain to the quantum Minkowski space
(see e.g. \cite{KUL1}) and in few other cases \cite{VLA,VAL}.

In the present paper we propose different extensions of this twist element.
The suggested construction implies the existence (in the universal
enveloping algebra to be deformed) of a subalgebra ${\bf L}$ with special
properties of multiplication. This is a solvable subalgebra with at least
four generators. All simple Lie algebras except $sl(2)$ contain such ${\bf L}$
and in any of them a deformation induced by twist of ${\bf L}$ can be
performed. In particular we study a jordanian deformation of ${\cal U}
(sl(N))$,
reaching a closed form of deformed compositions lacking in \cite{GER}. Using
the notion of factorizible twist \cite{RSTS} we prove that the element $%
{\cal F}\in {\cal U} (sl(N))^{\otimes 2}$,
\begin{equation}
\label{twist-sl(N)}{\cal F}=\exp \{2\xi \sum_{i=2}^{N-1}E_{1i}\otimes
E_{iN}e^{-\sigma }\}\exp \{H\otimes \sigma \},
\end{equation}
where $x=E_{1N}$, $H=E_{11}-E_{NN}$, $\sigma =\frac 12\ln (1+2\xi x)$,
satisfies the twist equation. Hence, it defines a triangular deformation of $%
{\cal U}(sl(N))$. In such Hopf algebras deformed by jordanian twist the
subset of Cartan generators $ \left\{ E_{ii} - E_{jj} \right\} $ with
$ i<j; \quad i,j \neq 1,N $ remains untouched. Hence there is a possibility
to perform additional multiparametric deformation using Reshetikhin twist
\cite{RES2}.
The main ingredients of the quantum group theory \cite
{DRIN} are constructed: the universal $R$-matrix, the dual Hopf algebra
(quantized function algebra on $SL(N)$), the universal ${\cal T}$-matrix
(canonical element) for the subalgebra which induces the twist of ${\cal U}%
(sl(N))$ and the self-duality of ${\bf L}$. Cohomological interpretation
of the interrelation between the
Drinfeld-Jimbo (or standard) quantum algebra ${\cal U}_q(sl(N))$ and the
jordanian (or non-standard) one ${\cal U}_\xi (sl(N))$ is discussed. The
real form and the corresponding quantum linear space are given. We present
also further generalization in which the subalgebra ${\bf L}$ is substituted
by a certain type of inhomogeneous Lie algebras.

The connection of the Drinfeld-Jimbo deformation \cite{DRIN,JIMB} with the
jordanian deformation was already pointed out in \cite{GER}. The similarity
transformation
of the classical matrix $r_{DJ}$ performed by the operator $\exp (\xi
adE_{1N})$ (with the highest root generator $E_{1N})$ turns $r_{DJ}$
into the sum $r_{DJ}+\xi r_j$ \cite{GER}. Hence,
\begin{eqnarray}r_j=-\xi \left( H_{1N}\wedge E_{1N}+2\sum_{k=2}^{N-1}E_{1k}\wedge
E_{kN}.\right) ,
\end{eqnarray}is a classical $r$-matrix too, which defines corresponding
deformation. A singular contraction of the quantum
Manin plane $xy=qyx$ of ${\cal U%
}_q(sl(2))$ with the mentioned above transformation in the fundamental
representation $M=1+\theta \rho (E_{1N})$, $\theta =\xi
(1-q)^{-1}$ results
in the jordanian plane $x^{\prime }y^{\prime }
=y^{\prime }x^{\prime }+\xi {%
y^{\prime }}^2$ of ${\cal U}_{\xi}(sl(2))$
\cite{OGIEV}. Later, this singular contraction in the
fundamental representation of $sl(3)$ and $sl(N)$ was used in many papers
(cf \cite{ALI,ABD} and references therein). 
Let us point out that in our formulas 
we do not refer to any particular representation of deformed algebras.

The paper is organized as follows. After reminding briefly the basic 
material on twisting of Hopf algerbas (Sec.2), we construct an extended 
jordanian twist $\cal F$ for four generator Lie algebra and apply 
it to twist the universal enveloping algebra ${\cal U}(sl(N))$ (Sec.3). 
The next Section contains cohomological explanation of the connection 
between the Drinfeld-Jimbo and jordanian quantization. The main 
objects of the theory of quantum groups are constructed in Sec.5. 
Further generalization of the extended jordanian twist to a special  
class of inhomogenious Lie algebras and possible research topics are 
given in Sec.6 and in the Conclusion. 

\section{Twisting of Hopf algebras}

\label{Sec2}

A Hopf algebra ${\cal A}(m,\Delta ,\epsilon,S)$ with
multiplication $m\colon {\cal A}\otimes {\cal A}\to {\cal A}$,
coproduct $\Delta \colon {\cal A}\to {\cal A}%
\otimes {\cal A}$, counit $\epsilon \colon {\cal A}\to C$,
and antipode $S\colon
{\cal A}\to {\cal A}$ (see definitions in Refs.\cite{DRIN,FAD,ChP})
can be transformed \cite{D2} with an invertible element ${\cal F}\in {\cal A}%
\otimes {\cal A}$, ${\cal F}=\sum f_i^{(1)}\otimes f_i^{(2)}$ into a twisted
one ${\cal A}_t(m,\Delta _t,\epsilon ,S_t)$.
This Hopf algebra ${\cal A}_t$ has the
same multiplication and counit maps but the twisted coproduct and antipode
$$
\Delta _t(a)={\cal F}\Delta (a){\cal F}^{-1},\quad S_t(a)=vS(a)v^{-1},\quad
v=\sum f_i^{(1)}S(f_i^{(2)}),\quad a\in {\cal A}.
$$
The twisting element has to satisfy the identities
\begin{eqnarray}(\epsilon \otimes  id)({\cal F}) = (id \otimes  \epsilon)({\cal F})=1,
\end{eqnarray}%
\begin{eqnarray}
{\cal F}_{12}(\Delta \otimes  id)({\cal F}) =
{\cal F}_{23}(id \otimes  \Delta)({\cal F}),
\label{TE}
\end{eqnarray}where the first one is just a normalizing condition and
follows from the second relation modulo a non-zero scalar factor.

A quasitriangular Hopf algebra ${\cal A}(m,\Delta , 
\epsilon,S,{\cal R})$ has
additionally an element ${\cal R}\in {\cal A}\otimes {\cal A}$ (a universal $%
R$-matrix) satisfying \cite{DRIN}
\begin{eqnarray}(\Delta \otimes  id)({\cal R})={\cal R}_{13}{\cal R}_{23}, \quad
(id \otimes  \Delta)({\cal R})={\cal R}_{13}{\cal R}_{12}.
\label{UR}
\end{eqnarray}The coproduct $\Delta $ and its opposite $\Delta ^{{\rm op}}$
are related by the similarity transformation (twisting) with ${\cal R}$
$$
\Delta ^{{\rm op}}(a)={\cal R}\Delta (a){\cal R}^{-1},\quad a\in {\cal A,}
$$
and in this case the relation (\ref{TE}) is just the Yang-Baxter equation.

A twisted quasitriangular quantum
algebra ${\cal A}_t(m,\Delta _t,\epsilon,S_t,%
{\cal R}_t)$ has the twisted universal $R$-matrix
\begin{eqnarray}{\cal R}_t=\tau  ({\cal F})\,{\cal R}\,{\cal F}^{-1},
\label{Rt}
\end{eqnarray}
where $\tau$ means permutation of the tensor factors: 
$\tau (f\otimes g)=(g\otimes f)\,, \tau ({\cal F})={\cal F}_{21}$.

Although, in principle, the possibility to quantize an arbitrary Lie
bialgebra has been proved \cite{ETI}, an explicit formulation of Hopf
operations remains a nontrivial task. In particular, the knowledge of
explicit form of the twisting cocycle is a rare case even for classical
universal enveloping algebras, despite of advanced Drinfeld's theory \cite
{D3}. Most of such explicitly known twisting elements have the factorization
property with respect to comultiplication (cf. (\ref{UR}))
$$
(\Delta \otimes id)({\cal F})={\cal F}_{23}{\cal F}_{13}\quad \mbox{or}\quad
(\Delta \otimes id)({\cal F})={\cal F}_{13}{\cal F}_{23}\,,
$$
and similar property involving $(id\otimes \Delta )$. To satisfy the twist
equation, these identities are combined with additional requirement 
${\cal F}_{12}{\cal F}_{23}={\cal F}_{23}{\cal F}_{12}$ 
or the Yang-Baxter equation on ${\cal F}$ \cite{RSTS,RES2}.

An important subclass of factorizable twists consists of elements
satisfying the following equations
\begin{equation}
\label{f-twist1}(\Delta \otimes id)({\cal F})={\cal F}_{13}{\cal F}_{23}\,,
\end{equation}
\begin{equation}
\label{f-twist2}(id\otimes \Delta _t)({\cal F})={\cal F}_{12}{\cal F}_{13}\,.
\end{equation}
It is easy to see that the universal $R$-matrix ${\cal R}$ satisfies these
equations, for $\Delta _t=\Delta ^{op}$. Another well developed case is the
jordanian twist of $sl(2)$  with ${\cal F}$ (\ref{og-twist}) \cite{OGIEV}. 
Due to the fact that the Cartan element $h$ is primitive in $sl(2):\Delta
(h)=h\otimes 1+1\otimes h$, and $\sigma  $ is primitive in the jordanian $%
{\cal U}_\xi (sl(2)):\Delta _t(\sigma   )=\sigma \otimes 1+1\otimes \sigma $,
one gets
$$
(\Delta \otimes id)e^{h\otimes \sigma }=e^{h\otimes 1\otimes \sigma
}e^{1\otimes h\otimes \sigma }\,,
$$
$$
(id\otimes \Delta _t)e^{h\otimes \sigma }=e^{h\otimes \sigma \otimes
1}e^{h\otimes 1\otimes \sigma }\,.
$$
It will be shown in the next Sec.3 that the element ${\cal F}$ (\ref
{twist-sl(N)}) also satisfies the factorization equations (\ref{f-twist1}),(%
\ref{f-twist2}) and can be used to twist the universal enveloping algebra of
$sl(N)$.

Let us mention, that the composition of appropriate twists can be defined
${\cal F} = {\cal F}_2 {\cal F}_1$. The element ${\cal F}_1$ has to
satisfy the twist equation with the coproduct of the original Hopf algebra,
while ${\cal F}_2$ must be its solution for $\Delta_{t_1}$ of the
intermediate Hopf algebra twisted by ${\cal F}_1$.
In particular, if ${\cal F}$ is a solution to the twist equation (\ref{TE})
then ${\cal F}^{-1}$ satisfies this equation with $\Delta \rightarrow \Delta
_t$.

\section{Factorizable twists}

Now we shall propose a new factorizable twist similar to (\ref{og-twist})
and defined on the abstract set of generators.

Let ${\bf L}$ be a four dimensional Lie algebra with generators $\left\{
H,A,B,E\right\} $ containing ${\bf B_{+}}$ and representable in a form of
semidirect sum of one dimensional space $V_H$ with basic element $H$ and a
Heisenberg subalgebra ${\cal H}(A,B,E)$ : ${\bf L=}V_H\vdash {\cal H}$ :
\begin{equation}
\label{l-comm}
\begin{array}{l}
\left[ H,E\right] =2E, \\
\left[ H,A\right] =\alpha A,\quad \left[ H,B\right] =\beta B,\qquad \alpha
+\beta =2, \\
\left[ E,A\right] =\left[ E,B\right] =0, \\
\left[ A,B\right] =\gamma E.
\end{array}
\end{equation}
Extending the twist deformation ${\cal U}_t\left( {\bf B_{+}}\right) $
performed by
$$
\Phi = \exp(\frac 12 H \otimes \ln(1+ \gamma E)) = e^{H \otimes \sigma} 
$$
to the universal enveloping ${\cal U}
\left( {\bf L}\right) $ one gets
the twisted algebra ${\cal U}_\Phi \left( {\bf L}\right) $ .
It retains the initial
multiplication defined by (\ref{l-comm}) while its coproduct $\Delta_{\Phi}
= \Phi \, \Delta \, \Phi^{-1}$ becomes noncocommutative:
\begin{equation}
\label{l-phi-cop}
\begin{array}{l}
\Delta _\Phi \left( H\right) =H\otimes e^{-2\sigma }+1\otimes H, \\
\Delta _\Phi \left( A\right) =A\otimes e^{\alpha \sigma }+1\otimes A, \\
\Delta _\Phi \left( B\right) =B\otimes e^{\beta \sigma }+1\otimes B, \\
\Delta _\Phi \left( E\right) =E\otimes e^{2\sigma }+1\otimes E,
\end{array}
\end{equation}
We shall show that the algebra ${\cal U}
\left( {\bf L}\right) $ allows a more
complicated twist deformation containing $\Phi $ as a factor.

\underline{Proposition.} The element
\begin{equation}
\label{twist}{\cal F}=\Phi \Phi _1=\exp (H\otimes \sigma )\exp (A\otimes
Be^{-2\sigma })
\end{equation}
is a twist for ${\cal U}  \left( {\bf L}\right) $ .

\underline{Proof.} We shall show that ${\cal F}=\Phi \Phi _1$ belongs to the
subclass defined by the equations (\ref{f-twist1},\ref{f-twist2}). The
equation (\ref{f-twist1}) is obviously true: $H$ and $A$ are the primitive
elements and $B$ commutes with $\sigma $ in ${\cal U}
\left( {\bf L}\right) $. To
check the second equation (\ref{f-twist2}) let us consider the coproducts 
$\Delta _{{\cal F}}\left( \sigma \right) $ 
and $\Delta _{{\cal F}}\left(B\right) $ . 
It is known that in twisted (by $\Phi $ ) universal enveloping 
of Borel subalgebra the element $\sigma $ is primitive \cite{OGIEV}. The
element $\sigma $ commutes not only with $B$ but also with $A$, 
so $\sigma $ remains primitive with respect to 
$\Delta _{{\cal F}}$. Using the properties 
of ''roots'' $\alpha -2= - \beta $ the twisted coproduct of $B$ can be
written in the following form
$$
\begin{array}{c}
\Delta _{
{\cal F}}\left( B\right) =\exp \left( {\rm ad}\left( A\otimes Be^{-\beta
\sigma }\right) \right) \circ \exp \left( {\rm ad}\left( H\otimes \sigma
\right) \right) \circ \left( B\otimes 1+1\otimes B\right) = \\ =\exp \left(
{\rm ad}\left( A\otimes Be^{-\beta \sigma }\right) \right) \circ \left(
B\otimes e^{\beta \sigma }+1\otimes B\right)  \\ =\exp \left( {\rm ad}\left(
A\otimes Be^{-\beta \sigma }\right) \right) \circ \left( B\otimes e^{\beta
\sigma }\right) +1\otimes B.
\end{array}
$$
From (\ref{l-comm}) one can see that 
$\left( {\rm ad}_A\right) ^2\circ B=0$. 
So the obtained expression can be simplified, 
$$
\Delta _{{\cal F}}\left( B\right) =B\otimes e^{\beta \sigma }+\left(
1+\left[ A,B\right] \right) \otimes B=B\otimes e^{\beta \sigma }+e^{2\sigma
}\otimes B.
$$
Now using the coproduct%
$$
\Delta _{{\cal F}}\left( Be^{-2\sigma }\right) =Be^{-2\sigma }\otimes
e^{-\alpha \sigma }+1\otimes Be^{-2\sigma }
$$
one can easily see that%
$$
\exp \left( {\rm ad}\left( H\otimes 1\otimes \sigma \right) \right) \circ
\left( A\otimes Be^{-2\sigma }\otimes e^{-\alpha \sigma }\right) =A\otimes
Be^{-2\sigma }\otimes 1.
$$
The latter garantes the validity of the equation (\ref{f-twist2}) for the
twisting element ${\cal F}$ . $\bullet $

The deformed algebra ${\cal U}_{{\cal F}}\left( {\bf L}\right) $ has initial
commutation relations generated by (\ref{l-comm}) and twisted coproducts:
\begin{equation}
\label{l-t-cop}
\begin{array}{l}
\Delta _{
{\cal F}}\left( H\right) =H\otimes e^{-2\sigma }+1\otimes H-2A\otimes
Be^{\left( \alpha -4\right) \sigma }, \\ \Delta _{
{\cal F}}\left( A\right) =A\otimes e^{-\beta \sigma }+1\otimes A, \\ \Delta
_{
{\cal F}}\left( B\right) =B\otimes e^{\beta \sigma }+e^{2\sigma }\otimes B,
\\ \Delta _{{\cal F}}\left( E\right) =E\otimes e^{2\sigma }+1\otimes E,
\end{array}
\end{equation}

Let us rewrite the twist element ${\cal F}$ in the reverse order:
\begin{equation}
\label{twist-rev}{\cal F}=\widetilde{\Phi _1}\Phi =\exp (A\otimes Be^{-\beta
\sigma })\exp (H\otimes \sigma )
\end{equation}
Now we know that both ${\cal F}$ and $\Phi $ are twists for ${\cal U}
\left( {\bf L}%
\right) $ and both satisfy the equations (\ref{f-twist1}),(\ref{f-twist2}).
Hence $\widetilde{\Phi _1}$ is also a twist element with respect to the
algebra ${\cal U}_\Phi \left( {\bf L}\right) $. 
Using the coalgebra relations (\ref{l-phi-cop}) it is easy to check 
that $\widetilde{\Phi _1}$ satisfies the
general twist equation (\ref{TE}), 
$$
\left( \widetilde{\Phi _1}\right) _{12}\left( \Delta _\Phi \otimes {\rm id}%
\right) \widetilde{\Phi _1}=\left( \widetilde{\Phi _1}\right) _{23}\left(
{\rm id}\otimes \Delta _\Phi \right) \widetilde{\Phi _1}.
$$
Note that contrary to the properties of ${\cal F}$ and $\Phi $ this twist ($
\widetilde{\Phi _1}$) does not belong to the subclass of factorizable twists
defined by the equations (\ref{f-twist1}),(\ref{f-twist2}).

Subalgebras of the type ${\bf L}$ exist in a large class of Lie algebras.
They can also be found in any simple Lie algebra of rank greater than 1.
Such simple algebras contain at least one pair of roots $\lambda _1$and $%
\lambda _2$ such that $\lambda _3=\lambda _1+\lambda _2$ is also a root. The
corresponding generators $X_1,X_2,X_3$ together with the Cartan element $H_3$
dual to the root $\lambda _3$ form the subalgebra equivalent to ${\bf L}$.
As we have shown above such subalgebra can be twisted with 
the element ${\cal F}$ and the corresponding 
deformation can be extended to the whole 
algebra ${\cal U}  $ and its twisted 
version ${\cal U}_{{\cal F}}$ can be thus constructed. 

We shall demonstrate the deformations generated by these twists in case of
simple algebras of series $A_{N-1}$. 
For our purposes it will be convenient to use 
the canonical basis of $gl(N)$ for the compositions of ${\cal U}(sl(N))$
\begin{equation}
\label{sln-comm}\,[E_{ik},E_{lm}]=\delta _{kl}E_{im}-\delta
_{im}E_{lk},\qquad i,k,l,m=1,...N.
\end{equation}
The Cartan elements of ${\cal U}  (sl(N))$
will be fixed as $H_{ik}=E_{ii}-E_{kk}$ .

Let $H \in {\bf L}$ be identified with the Cartan element 
dual to the highest 
root of $sl(N)$, this root will be denoted by $\lambda _H$. Collecting all
the pairs of roots with the property $\lambda _H=\lambda _1+\lambda _2$ one
can get the multiparametric twist of the type ${\cal F}$ with
$$
H=H_{1N},\ E=E_{1N},
$$
\begin{equation}
\label{aform}A=\sum_k\left( a^{1k}E_{1k}+a^{kN}E_{kN}\right) ,
\end{equation}
\begin{equation}
\label{bform}B=\sum_k\left( b^{1k}E_{1k}+b^{kN}E_{kN}\right) ,
\end{equation}
$$
\left\{ \ k=2,\dots ,N-1;\ a^{mn},b^{nm}\in {\sf C}\quad \right\} ,
$$
Here it is convenient to put $\gamma =2\xi $ ,%
$$
\sigma (E)=\frac 12\ln \left( 1+2\xi E\right) .
$$
In these terms the consistency condition would take the form
\begin{equation}
\label{a-b-corr}[A,B]=e^{2\sigma }-1 = 2\xi E 
\end{equation}
and the only nontrivial commutator of $\sigma$ with the basic elements of
${\bf L}$ is
\begin{equation}
\label{h-sigma}[H, \sigma]= 1- e^{-2\sigma}.
\end{equation}
According to the Proposition the element
\begin{equation}
\label{sln-twist}
{\cal F}=\Phi \Phi _1=\exp (H \otimes \sigma )\exp
(A\otimes Be^{-2\sigma })
\end{equation}
is a twist of ${\cal U}(sl(N))$. 
Using the particular properties of ${\bf L}$ one 
can apply the Cambell-Hausdorf formula to rewrite the twisting element in
the following form  
\begin{equation}
\label{sln-twist1}
\begin{array}{c}
{\cal F}=\exp (A\otimes B e^{-\sigma })\exp (H \otimes \sigma ) \\ =\exp
\left( H \otimes \sigma +A\otimes B \sigma e^{-2\sigma} 
\left(1- e^{-\sigma} \right)^{-1}\right) .
\end{array}
\end{equation}

\underline{Note}. Any number of factors of the type $\Phi _1$ can appear in
the expression (\ref{sln-twist}):
\begin{equation}
\label{fact-f}{\cal F}=\Phi \prod_j\Phi _j=\exp (H \otimes \sigma
)\prod_j\exp (A_j\otimes B_je^{-2\sigma })
\end{equation}
with $A_j$ and $B_j$ as in (\ref{bform}) and (\ref{aform}) and subject to
the additional conditions
$$
[A_{j_1},A_{j_2}]=[B_{j_1},B_{j_2}]=0,
$$
while the correlation equation (\ref{a-b-corr}) takes the form
\begin{equation}
\label{fact-corr}[A_j,B_k]=\delta _{jk}\left( e^{2\sigma }-1\right) . 
\, \bullet 
\end{equation}

Using the twist (\ref{sln-twist}) with the sole factor $\Phi _1$ one gets
the maximal number of free parameters -- the relation (\ref{a-b-corr})
imposes the only condition on the coefficients $a$'s and $b$'s,
\begin{equation}
\label{a-bcorr1}\sum_{k=2}^{N-1}\left( a^{1k}b^{kN}-a^{kN}b^{1k}\right)
=2\xi .
\end{equation}
On the contrary, supplying ${\cal F}$ with the maximal number ($N-2$) of
factors $\Phi _j$ one gets the $\left( N-2\right) ^2$ conditions (\ref
{fact-corr}). In particular one can satisfy $\left( N-2\right) \left(
N-3\right) $ of these conditions using the basic relations (\ref{sln-comm})
and the specific choice of $A_j$ and $B_j$ (one root $\lambda _j$ for each
factor $\Phi _j$ ):
\begin{equation}
\label{a-b-sln}
\begin{array}{c}
A_j=a^{1j}E_{1j}+a^{jN}E_{jN},\quad B_j=\left(
b^{1j}E_{1j}+b^{jN}E_{jN}\right), \\
\left( \mbox{with no summation on }j \, \right) . 
\end{array}
\end{equation}
Here the essential relations rest
\begin{equation}
\label{a-bcorr2}a^{1j}b^{jN}-a^{jN}b^{1j}=2\xi ,\quad \left\{
j=2,...N-1\right\} .
\end{equation}
Equation (\ref{a-bcorr1}) (as well as (\ref{a-bcorr2})) shows that it is
natural to renormalize the element $A$ (or the elements $A_j$ ) putting%
$$
A=2\xi \widetilde{A}
$$
so that%
$$
\ \left[ \widetilde{A},B\right] =E.
$$
In these notations the twisting elements
\begin{equation}
\label{sln-twist-tilde}{\cal F}=\exp (H \otimes \sigma )\exp (2\xi
\widetilde{A}\otimes Be^{-2\sigma }),
\end{equation}
\begin{equation}
\label{fact-f-tilde}{\cal F}=\exp (H \otimes \sigma )\prod_j\exp (2\xi
\widetilde{A_j}\otimes B_je^{-2\sigma })
\end{equation}
have the trivial limit $\lim_{\xi \rightarrow 0} {\cal F} = 1$. So does the
universal ${\cal R}$-matrix (${\cal R}={\cal F}_{21}{\cal F}^{-1}$) and one
can easily write down the corresponding classical $r$-matrices
\begin{equation}
\label{r-mat}r=-\left( H \wedge E\ +\ 2\widetilde{A}\wedge B\right)
\end{equation}
or
\begin{equation}
\label{r-mat-f}r=-\left( H \wedge E\ +\ 2\sum \widetilde{A}_j\wedge
B_j\right) .
\end{equation}
Their form clearly indicates that twisting by ${\cal F}$ corresponds to
the quantization of the self-dual Lie bialgebra $( {\bf L}, {\bf L}^*
\approx {\bf L} )$ just as in the case of the jordanian twist 
of ${\bf B}(1)$ \cite{OGIEV,VLA}. 
The same is true for the twisted Hopf algebra ${\cal U}_{{\cal F}}({\bf %
B^{\vee })}$, it is self-dual. We shall discuss this property in the next
Section and prove it in the Section 5 where the canonical element will
be constructed.

For the special case of ${\cal U}(sl(N))$ according to
the Proposition the 
following form of twisting element ${\cal F}$ can be chosen
$$
{\cal F}=\exp \left( H_{1N}\otimes \sigma \right)\prod_{j = 2}^{N-1}
\exp \left( 2\xi E_{1j}\otimes E_{jN}e^{-2\sigma }\right) .
$$
This twist of ${\cal U}(sl(N))$ is generated
by the twist of ${\cal U}  ({\bf L)}$ 
(here ${\bf L}$ is the restricted Borel subalgebra 
${\bf B^{\vee }}$ of $sl(N)$ 
with the basic elements $\left\{ H_{1N},E_{1N},E_{1j},
E_{jN} \right\}_{j=2, \dots , N-1} $ )
leading to the Hopf algebra ${\cal U}_{\xi}({\bf B}^{\vee })$ 
with the initial commutation relations (as in (\ref{sln-comm})), 
the twisted coproducts
\begin{equation}
\label{twisted-b}
\begin{array}{c}
\Delta _{
{\cal F}}H_{1N}=H_{1N}\otimes e^{-2\sigma }+1\otimes H_{1N}-4\xi
\sum_{j=2}^{N-1} E_{1j}\otimes E_{jN}e^{-3\sigma }, \\
\Delta _{{\cal F}}E_{1i}
=E_{1i}\otimes e^{-\sigma }+1\otimes E_{1i}, \\
\Delta _{{\cal F}}E_{iN}=
E_{iN}\otimes e^\sigma +e^{2\sigma }\otimes E_{iN}, \\
\Delta _{{\cal F}}E_{1N}=E_{1N}\otimes e^{2\sigma }+1\otimes E_{1N},
\end{array}
\end{equation}
antipodes
\begin{equation}
\label{twisted-Sb}
\begin{array}{c}
S_{\cal F}\left( \sigma \right)
=- \sigma ,\qquad S_{\cal F}\left( E_{1i}\right)
=-E_{1i}e^\sigma , \\
S_{\cal F}\left( E_{iN}\right)
=-E_{iN}e^{-3\sigma },\quad S_{\cal F}\left( E_{1N}\right)
=-E_{1N}e^{-2\sigma }, \\
S_{\cal F}\left( H_{1N}\right) =-H_{1N}e^{2\sigma }-4\xi
\sum_{j=2}^{N-1}E_{1j}E_{jN}
\end{array}
\end{equation}
and the universal ${\cal R}$-matrix of the form
\begin{equation}
\label{R-sl(3)}
\begin{array}{c}
{\cal R}={\cal F}_{21}{\cal F}^{-1} \\
=\prod_j \exp \left( 2\xi E_{jN}e^{-\sigma
}\otimes E_{1j}\right) \exp \left( \sigma \otimes H_{1N}\right) \exp \left(
-H_{1N}\otimes \sigma \right) \prod_{j} \exp \left( -2\xi E_{1j}\otimes
E_{jN}e^{-\sigma }\right).
\end{array}
\end{equation}
The coproducts and antipodes for other elements of ${\cal U}_{\xi}(sl(N))$
can be calculated using the standard formulas. The obtained expressions are
rather cumbersome. Thus, for example, in the case of ${\cal U}_{\xi}(sl(3))$
the coproduct of $E_{32}$ looks like
$$
\begin{array}{l}
\Delta _{
{\cal F}}E_{32}=E_{32}\otimes e^{-\sigma }+1\otimes E_{32}\\ +
\xi H_{13} \otimes E_{12} e^{-2\sigma} + 2\xi E_{12}
\otimes H_{23} e^{-\sigma} \\
- \xi H_{13}E_{12} \otimes (e^{-\sigma} - e^{-3\sigma})\\
-4\xi^2 E_{12} \otimes E_{23}E_{12} e^{-3\sigma}\\
-4\xi^2 E_{12}^2 \otimes E_{23} e^{-4\sigma} \,. 
\end{array}
$$
Twisting the coproducts is acting by the exponential of the adjoint operator
defined on the tesor product ${\cal U}(sl(N))
\otimes {\cal U}(sl(N))$. One can check that this operator is nilpotent and
all the twisted coproducts can be expressed
through the finite number of its powers.

\section{Connections between standard and jordanian deformations}

It is well known that some sorts of jordanian deformations can be treated as
limiting structures for certain sequences of standard quantizations \cite
{GER,OGIEV,ALI,ABD}. As will be shown below this is due to the specific
properties of Lie bialgebras involved in the quantizations. These properties
are more transparent when formulated for quantum groups rather than for
quantum algebras. For this reason in the current Section we use the dual
picture to treat Lie bialgebraic characteristics.

The generators of the standard (FRT-deformed) quantum group 
$Fun_h(SL(N))$ ($h=\ln q$) will be described by the entries of 
the $N\times N$-matrix $T$. Let $T$ be
subject to the similarity transformation with the matrix
\begin{equation}
\label{transM}M=1+\frac \xi {q-1}\rho \left( E_{1N}\right)
\end{equation}
(for the generators the canonical coproduct ( $\Delta T=T\stackrel{.}{%
\otimes }T$) is conserved). As far as $q\neq 1$ the transformed quantum
group $Fun_{h;\xi }(SL(N))$ is equivalent to the original one. Compare the
corresponding Lie bialgebras: $\left( g,g_{h;0}^{*}\right) =\left(
sl(N),\left( sl(N)\right) ^{*}\right) $ and $\left( g,g_{h;\xi }^{*}\right) $%
. Here the Lie algebra $g=sl(N)$ is not changed, the transformation $%
T\rightarrow MTM^{-1}$ does not touch the canonical coproduct for the
generators of the Hopf algebra $Fun_h(SL(N))$. Only the second Lie
multiplication $\left( \mu _{h;0}^{*}:V_{g^{*}}\wedge V_{g^{*}}\rightarrow
V_{g^{*}}\right) $ changes:%
$$
\mu _{h;0}^{*}\rightarrow \mu _{h;\xi }^{*}.
$$
The structure of the similarity transformation shows that the new Lie
product decomposes as:
\begin{equation}
\label{sum-mu'}\mu _{h;\xi }^{*}=\mu _{h;0}^{*}+\xi \mu ^{\prime }.
\end{equation}
The component $\mu ^{\prime }$ is fixed by the commutation relations that
can be extracted from the transformed $RTT=TTR$ equations. For this purpose
one has to change the coordinate functions of $SL(N)$ arranged in matrix $T$
for the exponential ones $T=\exp (\epsilon Y)$ and also change the
parameters $h\longmapsto \epsilon h,\ \xi \longmapsto \epsilon \xi $ .
Tending $\epsilon $ to zero one gets both summands in (\ref{sum-mu'}). The
second one of them looks as follows:
\begin{equation}
\label{dual-sln1}
\begin{array}{c}
\mu ^{\prime }\left( Y_{1k},Y_{ij}\right) =2\delta _{ik}Y_{Nj},
{\rm \ for}\ k,j<N;\quad i>1, \\ \mu ^{\prime }\left( Y_{ij},Y_{lN}\right)
=-2\delta _{jl}Y_{Nj},\
{\rm for}\ j<N;\quad i,l>1, \\ \mu ^{\prime }\left( Y_{ij},Y_{1N}\right)
=-\delta _{j1}Y_{i1}-\delta _{iN}Y_{Nj},\ {\rm for}\ j<N;\quad i>1,
\end{array}
\end{equation}
\begin{equation}
\label{dual-sln2}
\begin{array}{c}
\mu ^{\prime }\left( Y_{1i},Y_{1N}\right) =-Y_{1i},\
{\rm for}\ i>1, \\ \mu ^{\prime }\left( Y_{1N},Y_{kN}\right) =Y_{kN},\
{\rm for}\ k<N, \\ \mu ^{\prime }\left( Y_{11},Y_{1N}\right) =\mu ^{\prime
}\left( Y_{1N},Y_{NN}\right) =-\left( Y_{11}-Y_{NN}\right) , \\
\mu ^{\prime }\left( Y_{1i},Y_{1k}\right) =\delta _{i1}Y_{Nk},
{\rm \ for}\ k,i<N, \\ \mu ^{\prime }\left( Y_{iN},Y_{kN}\right) =-\delta
_{kN}Y_{i1},
{\rm \ for}\ k,i>1, \\ \mu ^{\prime }\left( Y_{1i},Y_{kN}\right) =\delta
_{i1}Y_{k1}-\delta _{kN}Y_{Ni}-2\delta _{ik}\left( Y_{11}-Y_{NN}\right) ,\
{\rm for}\ i<N;\quad k>1.
\end{array}
\end{equation}
Here for simplicity of exposition we use the canonical $gl(N)$-basis. One
can check that this deforming function $\mu ^{\prime }$ not only defines the
infinitesimal deformation of $\mu _{h;0}^{*}$ but is itself a Lie
multiplication.

Consider the decomposition (\ref{sum-mu'}) as a deformation equation for the
original dual Lie algebra $g_{h;0}^{*}$. Its main property is that $\mu
^{\prime }$ does not depend on $h$ or $\xi $. So the transformed law has the
form
\begin{equation}
\label{sum-mu}\mu _{h;\xi }^{*}=\mu _{h;0}^{*}+\mu _{0;\xi }^{*}.
\end{equation}
This means that $\mu _{h;\xi }^{*}$ is a Lie multiplication deformed in the
first order. Both summands are Lie maps and at the same time can be
considered as deforming functions of each other. As a result both deforming
functions are Hochschild 2-cocycles for the corresponding Lie algebras ($%
g_{0;\xi }^{*}$ with the multiplication $\mu _{0;\xi }^{*}$ and $g_{h;0}^{*}$
defined by $\mu _{h;0}^{*}$)
$$
\begin{array}{c}
\mu _{h;0}^{*}\in Z^2\left( g_{0;\xi }^{*},g_{0;\xi }^{*}\right) , \\
\mu _{0;\xi }^{*}\in Z^2\left( g_{h;0}^{*},g_{h;0}^{*}\right) .
\end{array}
$$
The equivalence of the algebraic structures in $Fun_{h;\xi }(SL(N))$ and $%
Fun_h(SL(N))$ (for $h\neq 0$) signifies that $\mu _{0;\xi }^{*}$ is in fact
a coboundary,%
$$
\mu _{0;\xi }^{*}\in B^2\left( g_{h;0}^{*},g_{h;0}^{*}\right) .
$$
On the contrary, the composition $\mu _{h;0}^{*}$ corresponds to a
nontrivial cohomology class%
$$
\mu _{h;0}^{*}\in H^2\left( g_{0;\xi }^{*},g_{0;\xi }^{*}\right) ,
$$
the deformation of $\mu _{0;\xi }^{*}$ by $\mu _{h;0}^{*}$ is essential
\cite{NIJ}.

Notice that the multiplication maps here have certain cohomological
properties also with respect to cochain complex $C$ of maps  $%
C^n:\bigwedge^nV_g\rightarrow V_g\wedge V_g$, where the $g$-module is chosen
to be $\bigwedge^2V_g$ with the canonically extended adjoint action on it.
The dualization of spaces $V_g\Leftrightarrow V_{g^{*}}$ converts the map $%
\mu ^{*}$ into the chain $\mu ^{*}\in C^1\left( g,g\wedge g\right) $ . As it
was mentioned above the initial coproduct for the generators of $Fun_h(SL(N))
$ rests unchanged under the transformation. All the Lie algebras $g_{h;\xi
}^{*}$ are dual to one and the same $g=sl(N)$. Thus both $\mu _{h;0}^{*}$ and
$\mu _{0;\xi }^{*}$ are 1-cocycles for the complex $C$ .

This set of characteristics necessarily indicates that the classical 
$r$-matrix of \\
 ${\cal U}_{h;\xi }\left( sl(N)\right) \approx $ 
$\left( Fun_{h;\xi}(SL(N))\right)^{*} $ 
must also exhibit this decomposition property: 
\begin{equation}
\label{r-sln}
\begin{array}{c}
r_{h;\xi }=r_{h;0}+r_{0;\xi }=
\frac hN\left( \sum_{k=1}^{N-1}k\left( N-k\right)
H_{k,k+1}\otimes H_{k,k+1}\right.  \\
+\left. \sum_{k<l}\left( N-l\right) k\left( H_{k,k+1}\otimes
H_{l,l+1}+H_{l,l+1}\otimes H_{k,k+1}\right) \right)  \\
+2h\sum_{k<l}\left( E_{lk}\otimes E_{kl}\right)  \\
-\xi H_{1N}\wedge E_{1N}-2\xi \sum_{k=2}^{N-1}E_{1k}\wedge E_{kN}.
\end{array}
\end{equation}
In the limit $h\rightarrow 0$ one gets the element 
\begin{eqnarray}
\lim _{h\rightarrow 0}r_{h;\xi }=r_{0;\xi }=-\xi
\left( H_{1N}\wedge
E_{1N}+2\sum_{k=2}^{N-1}E_{1k}\wedge E_{kN}.\right) ,
\label{CYBE0}
\end{eqnarray}
that coincides with $r$-matrix that can be obtained from ${\cal R}$
presented above (\ref{R-sl(3)}). So the jordanian quantum group 
$Fun_{0;\xi}(SL(N))$ has the same $R$-matrix as the twisted algebra 
${\cal U}_{\cal F}(sl(N))$ (with ${\cal F}$ 
as in (\ref{fact-f-tilde}) and 
$\widetilde{A_j}, B_j$ as in (\ref{a-b-sln}),(\ref{a-bcorr2})).

The $r$-matrices (\ref{r-sln}) and (\ref{CYBE0}) are known for a long time.
In \cite{GER} $r_{0;\xi }$ was obtained by applying ${\mbox{\rm ad}}%
_{E_{1N}}^{}$ to the canonical antisymmetric $r_{\wedge
}=\sum_{i<j}E_{ij}\wedge E_{ji}$. It was stressed that $r_{0;\xi }$ lay in
the boundary of the dense set of orbits of standard quantizations induced by
$r_{\wedge }$. The $r$-matrix (\ref{CYBE0}) was also obtained in the
discussion of conformal algebra deformations \cite{LUM}.

The $r$-matrix  $r_{0;\xi }$ is the element of the space ${\bf B}^{\vee
}\bigwedge {\bf B}^{\vee }.$ Its structure suggests the renumeration of the
basic elements of ${\bf B}^{\vee }$; we shall describe the
corresponding basis as the set
$$
\left\{ P_\alpha ,X_\beta \right\} _{\alpha ,\beta =1,\ldots ,N-1}\quad {\rm %
with}\quad \left\{
\begin{array}{c}
P_1=E_{1N},\ P_i=E_{iN}; \\
X_1=H_{1N},\ X_j=2E_{1j};
\end{array}
\right. \ i,j=2,\ldots ,N-1.
$$
In these notations $r_{0;\xi }$ takes the form%
$$
r_{0;\xi }=-\xi X_\alpha \wedge P_\alpha .
$$
The basic exponential coordinate functions $\left\{
Y_{1N},Y_{iN},Y_{11}-Y_{NN},Y_{1i}\right\} $ are chosen so that they are
canonically dual to those of $\left\{ P_\alpha ,X_\beta \right\} $ . Let us
apply the homomorphism
\begin{equation}
\label{r-hom}r_{0;\xi }:Y\rightarrow -\xi X_\alpha \wedge \left\langle
P_\alpha ,Y\right\rangle
\end{equation}
to the Lie algebra $\left( {\bf B}^{\vee }\right) ^{*}$ described by the last
six compositions $\mu ^{\prime }$ (see (\ref{dual-sln2})). As a result we
shall get the Lie algebra ${\bf B}^{\vee }$ . The significant fact is that (%
\ref{r-hom}) is an isomorphism, that is ${\bf B}^{\vee }\approx \left(
{\bf B}^{\vee }\right) ^{*}$. The twist ${\cal F}$ induces the self-dual Lie
bialgebra $\left( {\bf B}^{\vee },{\bf B}^{\vee }\right) $.

It is useful to compare this situation with that of a classical double of
dual Lie algebras $\left( g,g^{*}\right) $ .There the composition law of the
double can also be presented as a sum of two multiplications with
independent linear parameters. But in that case both summands are
cohomologically nontrivial. What is more important -- such parametrization
(and subdivision) can not be performed in only one algebra of a Lie
bialgebra $\left( g,g^{*}\right) $ corresponding to a classical double. In
fact these are the Lie bialgebras that can be parametrized in that case so
that their arbitrary linear combination is again a Lie bialgebra \cite{Lya}.
When a Lie bialgebra is nontrivialy decomposed (that is the decomposition
goes parallel in both dual algebras) the $r$-matrix for a linear
combination of Lie bialgebras doesn't inherit the decomposition property.

To clarify the contraction properties of 
$Fun_{h;\xi }(SL(N))$ let us 
consider the 1-parameter subvariety 
$\left\{ g_{h;1-h}^{*}\right\} $ of Lie 
algebras $g^{*}{}_{h;\xi }^{}$ (putting $\xi =1-h$ in (\ref{sum-mu})). 
Each dual pair $\left( sl(N),g_{h;1-h}^{*}\right) $ 
is a Lie bialgebra and thus 
is quantizable \cite{ETI}. The result is the set ${\cal A}_{s;h}$ of
deformation quantizations parametrized by $h$ 
and the deformation parameter $s$. 
This set can be considered smooth in the sense compatible with the
formal series topology \cite{BUR} -- close Lie bialgebras give rise to close
deformations. The 1-dimensional boundaries 
${\cal A}_{0;h}$ and ${\cal A}_{s;0}$ of ${\cal A}_{s;h}$ 
are formed respectively by the quantizations of 
$\left( sl(N),g_{1;0}^{*}\right) $ (the standard Lie bialgebra) and 
$\left(sl(N),g_{0;1}^{*}\right) $ 
(the jordanian one). Each internal point in ${\cal A}_{s;h}$ 
can be connected with a boundary by a smooth parametric
curve $a(u)$. One can choose it so that it starts in ${\cal A}_{0;h}$ and
ends in ${\cal A}_{s;0}$ . So a jordanian Hopf algebra obtained by twisting
deformation can be also treated as a limit point of a smooth 1-dimensional
subvariety $a(u)$ . This does not necessarily mean that this limit is a
faithful contraction --it may be impossible to attribute the curve $a(u)$ to
an orbit of some Hopf algebra in ${\cal A.}$ This is just what happens when
the transformation $M$ is applied to $Fun_{h;0}(SL(N))$. For every positive $%
h$ fixed the subset $\left\{ Fun_{h;\xi }(SL(N))\right\} $ is in the $GL(N^2)
$-orbit of the corresponding $Fun_{h;0}(SL(N))$. To attain the points $%
Fun_{0;\xi }(SL(N))$ one must tend $h$ to zero and this can be done only by
crossing the set of orbits refering to inequivalent Hopf algebras. These
specific interrelations of different types of quantizations where noted in
\cite{GER}. It was demonstrated for the case of $sl(N)$ that the standard
deformation $Fun_{h;0}(SL(N))$ can be accompanied by a smooth transformation
of a jordanian deformation so that the latter reaches the orbit of $%
Fun_{h;0}(SL(N))$ . Applying the operator $M$ to an element of the set $%
\left\{ Fun_{h;0}(SL(N))\right\} $ one gets an intersection point of an
orbit and of a curve parametrized by $\xi $ .

One of the principle conclusions is that the possibility to obtain the
jordanian deformation $Fun_{0;\xi }(SL(N))$ as a limiting transformation of
the standard quantum group --$Fun_{h;0}(SL(N))$ (and on the dual list to get
the twisted $q$-algebra ${\cal U}_{{\cal F}}(sl(N))$ as a limit of the variety of
standardly quantized algebras ${\cal U}_q(sl(N)))$ is provided by the fact that the
1-cocycle $\mu _{0;\xi }^{*}\in Z^1(sl(N),sl(N)\wedge sl(N))$ (that
characterizes the Lie bialgebra for ${\cal U}_{{\cal F}}(sl(N))$ ) is at the same
time the 2-coboundary $\mu _{0;\xi }^{*}\in B^2\left(
g_{h;0}^{*},g_{h;0}^{*}\right) $ the Lie algebra $g_{h;0}^{*}$ being the
standard dual of $sl(N)).$

\section{Canonical element and jordanian quantum space}

The set $\left\{ P_\alpha ,X_\beta \right\} $ forms the basis appropriate to
deal with the Lie bialgebras $\left({\bf L}, {\bf L}^{*}\right)$. 
To study the properties of $R$-matrix ${\cal R}$ and 
the canonical element ${\cal T}$ it is 
reasonable to perform the corresponding rearrangement of basis for the whole
Hopf algebra ${\cal U}_\xi \left({\bf B}^{\vee }\right) $. 
We shall consider the set 
\be
\label{newbase}
\left\{ z_k\right\} _{k=1,\ldots ,2\left( N-1\right) }=\left\{ x_\alpha ,\pi
_\beta \right\} _{\alpha ,\beta =1,\ldots ,N-1}
\ee
as the generators of ${\cal U}_\xi \left({\bf B}^{\vee }\right) $ with%
$$
\begin{array}{l}
x_1=H_{1N},\quad x_i=2E_{1i}, \\
\pi _1=\frac 1\xi \sigma =\frac 1{2\xi }\ln \left( 1+2\xi E_{1N}\right)
,\quad \pi _i=E_{iN}e^{-2\sigma }.
\end{array}
$$
The basis will be formed by
the set of ordered monomials:
\begin{equation}
\label{arrow-b}\left\{ z_{\vec k}\right\} _{{\vec k=\left\{ \vec m,\vec
n\right\} =\left\{ m_1,\ldots ,m_{N-1},n_1,\ldots ,n_{N-1}\right\} }%
}=\left\{ x_1^{m_1} \ldots x_{N-1}^{m_{N-1}} \pi _1^{n_1} \ldots 
\pi_{N-1}^{n_{N-1}}\right\} .
\end{equation}
In these terms the ${\cal R}$-matrix (\ref{R-sl(3)}) can be rewritten as
\begin{equation}
\label{ur-sln}{\cal R}=\prod_{\alpha =1,\ldots ,N-1}^{<}\exp (\pi _\alpha
\otimes \xi x_\alpha )\prod_{\alpha =1,\ldots ,N-1}^{>}\exp (-\xi x_\alpha
\otimes \pi _\alpha ),
\end{equation}
We shall use the standard Hopf algebra homomorphism ${\cal R}\colon {\cal A}%
^{*}\to {\cal A}_{-}$where in our case ${\cal A}$ is the twisted algebra $%
{\cal U}_\xi \left({\bf B}^{\vee }\right) $ and ''$-$'' indicates the opposite
multiplication. It would be appropriate to consider ${\cal R}$ as belonging
to ${\cal A}_{-}\otimes {\cal A}$ with the decomposition
\begin{equation}
\label{r-decomp}{\cal R}=\sum R^{{\vec k\vec l}}y_{{\vec k}}\otimes z_{{\vec
l}}.
\end{equation}
It is implied that the basic monomials $y_{{\vec k}}\in {\cal A  }_{-}$
contain the same sequences of generators $z_k$ as the corresponding basic
elements $z_{{\vec k}}\in {\cal A  }$ (see (\ref{arrow-b})) but the
multiplication that glue them is opposite to that of ${\cal A  }$. Let $%
\left\{ z^{_{\vec k}}\right\} =\left\{ x^{_\alpha },\pi ^{_\beta }\right\} $
be the canonical dual basis of ${\cal A  }^{*}$ ,$\left\langle z^{_{\vec
k}},z_{{\vec l}}\right\rangle =\delta _{{\vec l}}^{_{{\vec k}}}$. The
morphism ${\cal R}$ can be defined by its values on the basic elements:
\begin{equation}
\label{r-map}{\cal R}\left( z^{_{\vec k}}\right) =\sum R^{{\vec l\vec k}}y_{{%
\vec l}}
\end{equation}
Let us extract the first terms of the decomposition (\ref{r-decomp}) for the
${\cal R}$-matrix (\ref{ur-sln})
\begin{equation}
\label{r-expand}{\cal R}=1\otimes 1+R^{kl}z_k\otimes z_l+\ldots
\end{equation}
(Note that in such a presentation the second term is not proportional to the
classical $r$-matrix, the generators $z_l$ do not form a Lie algebra.) The
terms written explicitly in (\ref{r-expand}) are the only ones containing
the first powers of generators. Thus the images ${\cal R}\left( z^{_{\vec
k}}\right) $ are the linear combinations of the generators $z_{{\vec l}}$.
In our case the matrix $\left\{ R^{kl}\right\} $ is invertable,%
$$
R=\xi \left(
\begin{array}{cc}
0 & -I \\
I & 0
\end{array}
\right) \Rightarrow \left\{
\begin{array}{c}
{\cal R}\left( x^\alpha \right) =\xi \pi _\alpha  \\ {\cal R}\left( \pi
^\alpha \right) =-\xi x_\alpha
\end{array}
\right. ,
$$
$$
R^{-1}=-\frac 1{\xi ^2}R\quad \Rightarrow \left\{
\begin{array}{c}
{\cal R}^{-1}\left( x_\alpha \right) =-\frac 1\xi \pi _\alpha  \\ {\cal %
R}^{-1}\left( \pi _\alpha \right) =\frac 1\xi x^\alpha
\end{array}
\right. .
$$
Reversing the formula (\ref{r-map}) we get 
the expression for the elements of 
the dual basis in terms of generators $z^k_{\rule{2mm}{0mm}
|k=1,\ldots ,N-1}$:
\begin{equation}
\label{dual-base}z^{_{\vec k}}=\sum R^{{\vec l\vec k}}{\cal R}^{-1}\left( y_{%
{\vec l}}\right) =\sum R^{{\vec l\vec k}}\left( \left( R^{-1}\right)
_{k_11}z^{k_1}\right) ^{l_1}\ldots \left( \left( R^{-1}\right) _{k_{2\left(
N-1\right) }2\left( N-1\right) }z^{k_{2\left( N-1\right) }}\right)
^{l_{2\left( N-1\right) }}.
\end{equation}
The basic decomposition for the ${\cal R}$-matrix (\ref{ur-sln})
can be written explicitly,
\begin{equation}
\label{r-dec-base}
\begin{array}{c}
\begin{array}{rcl}
{\cal R} & = & \sum
\frac{(-\xi )^{|\vec n|} \xi^{|\vec m|}}{\vec m!\vec n!}x_1^{n_1}\ldots
x_{N-1}^{n_{N-1}}(\pi _1)^{m_1}\ldots (\pi _{N-1})^{m_{N-1}}\otimes  \\  &
\otimes  & x_1^{m_1}\ldots x_{N-1}^{m_{N-1}}(\pi _1)^{n_1}\ldots (\pi
_{N-1})^{n_{N-1}},
\end{array}
\\
\left| \vec n\right| =n_1+\ldots +n_{N-1};\quad \vec n!=n_1!n_2!\ldots
n_{N-1}!
\end{array}
\end{equation}
Here we used the inclusion ${\cal R}\in {\cal A  }_{-}\otimes {\cal A  }$ and
the fact that all the generators $\pi _\alpha $ commute. The structure of $%
{\cal R}$-morphism is clearly seen here. It states the one-to-one
correspondence between the basic monoms of ${\cal A  }^*$ and ${\cal A }_-$.
This evidently signifies that the
Hopf algebras ${\cal A  }^{*}$ and ${\cal A  }_{-}$ are equivalent.
One must also take into account that in our case ${\cal A  }_{-}$ is the
twisted universal enveloping algebra ${\cal U}_\xi
\left({\bf B}^{\vee }\right)$ with the opposite product.
Hence it is isomorphic to its multiplicatively inverse. The result
is
$$
{\cal A  }^{*} \approx {\cal A  }_{-}
\approx {\cal A  }
$$
The Hopf algebra $ {\cal U}_\xi \left({\bf B}^{\vee
}\right) $ is self-dual.

The structure constants $R^{{\vec l\vec k}}$ presented
in the decomposition (\ref{r-dec-base}) can be substituted in
the expression (\ref{dual-base}) to fix
explicitly the form of the dual basis.
Hence the canonical element ${\cal T}$ is completely defined %
$$
\begin{array}{rclc}
{\cal T} & =\sum_{{\vec k=\left( \vec m,\vec n\right) }}z^{_{{\vec k}%
}}\otimes z_{{\vec k}} & = & \sum \frac 1{\vec m!\vec n!}\left( \pi
^1\right) _{}^{n_1}\ldots \left( \pi ^{N-1}\right)
_{}^{n_{N-1}}(x^1)^{m_1}\ldots (x^{N-1})^{m_{N-1}}\ \otimes  \\
&  &  & \otimes \ x_1^{m_1}\ldots x_{N-1}^{m_{N-1}}\pi _1{}^{n_1}\ldots \pi
_{N-1}{}^{n_{N-1}}.
\end{array}
$$
We can recollect this expansion into the ordered product using the following
property of the ${\cal T}$-matrix: $(id\otimes S)({\cal T})={\cal T}^{-1}$.%
$$
\begin{array}{rlc}
{\cal T}^{-1} & = & \sum_{
{\left( \vec m,\vec n\right) }}\frac 1{\vec m!\vec n!}\left( \pi ^1\right)
_{}^{n_1}\ldots \left( \pi ^{N-1}\right) _{}^{n_{N-1}}(x^1)^{m_1}\ldots
(x^{N-1})^{m_{N-1}}\ \otimes  \\  &  & \otimes \ \left( S\left( \pi
_{N-1}\right) \right) ^{n_{N-1}}\ldots \left( S\left( \pi _1\right) \right)
^{n_1}\left( S\left( x_1\right) \right) ^{m_1}\ldots \left( S\left(
x_{N-1}\right) \right) ^{m_{N-1}}
\end{array}
$$
The antipodes used here can be easily found using the expressions
(\ref{twisted-Sb}) given in Sec.3: 
$$
\begin{array}{c}
S(\pi _1)=-\pi _1,\quad S(\pi _i)=-\pi _ie^{\xi \pi _1}, \\
S(x_1)=-x_1e^{2\xi \pi _1}-4\xi \sum_{}x_i\pi _ie^{2\xi \pi _1},\quad
S(x_i)=-x_ie^{\xi \pi _1}.
\end{array}
$$
Note that the homomorphic image in ${\cal A}^{*}$ of the abelian subalgebra
generated by elements $\left\{ \pi _\alpha \right\} $ is itself a
commutative subalgebra. This enables us to write the final expression for
the canonical element
\begin{equation}
\label{T-exp}{\cal T}=\prod_{{\alpha }}^{<}\exp (-x^\alpha \otimes
S(x_\alpha ))\prod_{{\alpha }}^{>}\exp (-\pi ^\alpha \otimes S(\pi _\alpha
)).
\end{equation}

The corresponding constructions for jordanian deformations of the Lie
superalgebra $sl(M|N)$ can be found in \cite{KUL3}.

Let us present a real form for ${\cal U}_{\cal F}(sl(N))$. 
We focus first on the subalgebra ${\bf B}^\vee $
in the general setting of the previous section and with the basis
$\{ z_k \}$ (see (\ref{newbase})).  The
anti-algebraic anti-linear transformation given on the generators by
$$\theta(x_\alpha) = - x_\alpha ,\quad \theta(\pi_\alpha) = \pi_\alpha
$$
respects the classical comultiplication and defines a real form on
${\cal U}({\bf L})$. 
At the same time, the twisting element ${\cal F}$ turns
into ${\cal F}^{-1}$. Henceforth, $\theta$ is a real form (cohomomorphic
and anti-homomorphic) for the twisted algebra 
${\cal U}_{\cal F}({\bf L})$ as well. 
Turning to the specific case of $sl(N)$, the question is whether
$\theta$ can be extended from the subalgebra ${\bf B}^\vee $ to
the entire $sl(N)$. This is possible, and the corresponding transformation
is
$$\theta(E_{ij})=-E_{ij},\quad i,j<N \quad\hbox{or}\quad i,j=N;\quad
  \theta(E_{kN})= E_{kN} ,\quad
  \theta(E_{Nk})= E_{Nk} ,\quad k <N.
$$
It is evident that $\theta$ is a Lie algebra anti-automorphism. The real 
form for $N = 2$ case of the jordanian 
${\cal U}_{\xi}(sl(2))$ was given in \cite{WOR}. 

Twisting of a symmetry Hopf algebra ${\cal A}$
of a manifold ${\cal M}$ induces deformation
of its whole geometry, so that the notion of symmetry is conserved
in the framework of the non-commutative geometry. Such deformation
includes that of function algebras ( vector bundles, $*$-structure, and
so on) expressing new  objects in terms of the untwisted ones
by explicit formulas involving twisting 2-cocycle ${\cal F}$ .
Here we present, as an application of the developed
jordanian-type quantization of $sl(N)$, the corresponding noncommutative
 space ${\cal M}_{\cal F}$. We deduce commutation relations
for generators of ${\cal M}_{\cal F}$,
and the differential calculus.
The basic formula connecting multiplications in 
${\cal A}$-modules ${\cal M}_{\cal F}$
and ${\cal M}$ (the twisted and the untwisted ones) is \cite{D3} 
\be
f*g={\cal F}^{-1}_{(1)}(f)\cdot{\cal F}^{-1}_{(2)}(g),\quad f,g\in {\cal M}.
\label{TM}
\ee
The star stands for the new product on ${\cal M}_{\cal F}$
defined through the old one '$\cdot$' and the element ${\cal F}$.
If ${\cal M}$ is classical, the twisting cocycle is represented by
a bidifferential operator according to the correspondent representation
of ${\cal F}$.
Thus ${\cal M}_{\cal F}$ and ${\cal M}$ coincide as linear spaces but
they are endowed with different algebraic structures. The transformation
is performed in such a way that the symmetry property
$h(f\cdot g)=h_{(1)}(f)\cdot h_{(2)}(g),
\quad h \in {\cal A},\quad f,g\in {\cal M}$,
is inherited by the twisted algebra ${\cal A}_{\cal F}$.

Let $x^\al, \> \al=1,\ldots N$, be the generators of
${\cal M}_{\cal F}$. To evaluate commutation relations among
the generators, it is sufficient to retain only the following terms:
$$
{\cal F}=1\oo 1 + \xi(x^1 \partial_1 - x^N \partial_N)\oo x^1 \partial_N
  +2\xi\sum_{k=2}^{N-1}  x^1 \partial_k \oo x^k \partial_N +\ldots ,
$$
with the rest of the series vanishing.
Resolving formula (\ref{TM}) (twisting is an invertible operation)
we come to
$$
{\cal F}_{(1)}(x^\mu)*{\cal F}_{(2)}(x^\nu)
=x^\mu \cdot x^\nu=x^\nu \cdot x^\mu
={\cal F}_{(1)}(x^\nu)*{\cal F}_{(2)}(x^\mu).
$$ 
This gives (commutators are understood in terms of 
the twisted product '$*$') \cite{GER} 
$$
\begin{array}{llllll}
[x^1,x^N]&=&\xi x^N* x^N, &
[x^i,x^k]&=&0       ,  \\[6pt]
[x^1,x^k]&=&2\xi x^k* x^N, &
[x^k,x^N]&=&0       .  \\[6pt]
\end{array}
$$
Hereafter (in this Section) the small Latin indices run from $2$ to $N-1$.
Similarly, for the contravariant entities $p_\mu$ we have
$$
\begin{array}{llllll}
[p_1,p_N]&=& \xi  p_1 * p_1, &
[p_i,p_k]&=&0       ,  \\[6pt]
[p_k,p_N]&=& 2\xi p_1* p_k, &
[p_1,p_k]&=&0       .  \\[6pt]
\end{array}
$$
Let us note that after the quantization the bases $\{ p_\mu \}$
and $\{ x^\mu \}$
are no longer conjugate. The invariant canonical element
turns out to be $x^\mu\cdot p_\mu= x^\mu* p_\mu+\xi x^N* p_1$.
Non-trivial cross-relations between coordinates and momenta are
expressed by
$$ 
[p_N,x^1]=\xi (p_N*x^N +2 \sum_{k=2}^{N-1} p_k*x^k + p_1*x^1 
+ \xi p_1*x^N) 
$$ 
$$ 
[p_1,x^1]=-\xi p_1* x^N, \quad [p_k,x^k]=-2 \xi p_1*x^N, 
\quad [p_N,x^N]=-\xi p_1*x^N . 
$$ 
Partial derivatives $\partial_\mu$ satisfies the same identities as
$p_\mu$, whereas the cross-relations are modified accordingly:
$$ 
[\partial_N,x^1]=\xi (x^N *\partial_N + 
2\sum_{k=2}^{N-1} x^k *\partial_k 
+ x^1* \partial_1 + \xi x^N*\partial_1) 
$$ 
$$ 
[\partial_1,x^1]=1 - \xi x^N*\partial_1\,, 
\quad [\partial_k,x^k]=1+2\xi  x^N*\partial_1\,, 
\quad [\partial_N, x^N]=1 - \xi x^N*\partial^N .
$$

\section{Group cocycles and twisting}

To generalize the construction of Sec.3
let us arrange the generators of ${\bf B}^\vee$
into the two sets $(H,A_j)$ and $(E,B_j)$ spanning two
mutually complement Lie sublagebras. We denote them  ${\bf H}$ and
${\bf H}^*$,
respectively, regarding as dual linear spaces.
Subalgebra ${\bf H}$ acts on ${\bf H}^*$, thus endowing
${\bf B}^\vee$ with the semidirect sum
${\bf L}={\bf H}\triangleright{\bf H}^*$ structure.
In this section we establish the cohomological properties of
the previous constructions in terms of
the Lie algebra ${\bf H}$ and its formal
Lie group ${\bf G}= \exp{\bf H}$.

Let $H_\mu $ be basic elements of a Lie algebra ${\bf H}$ and
$X^\nu$ be their conjugate.
Commutation relations in ${\bf H}$ are specified by the structure
constants $C_{\mu\nu}^\si$:
\be
[H_\mu,H_\nu]=C_{\mu\nu}^\si H_\si.
\label{HHrel}
\ee
Suppose a left action of  ${\bf H}$ on ${\bf H}^*$
\be
[H_\mu,X^\nu]=-L_{\mu\si}^\nu X^\si
\label{HXrel}
\ee
which enables us to build the semidirect sum
${\bf L}={\bf H}\triangleright{\bf H}^*$ where ${\bf H}^*$ is assumed
to be an abelian subalgebra.
The element
\be
r=X^\nu\otimes  H_\nu - H_\nu\otimes  X^\nu  \in
{\bf L}\wedge {\bf L}
\label{CYBE}
\ee
is a solution to the classical Yang-Baxter equation if and only if
\be
C_{\mu\nu}^\si = L_{\mu\nu}^\si-L_{\nu\mu}^\si.
\label{e1}
\ee
The structure constants $L_{\mu\nu}^\si$ defines also a left action
of ${\bf H}$ on itself according to the rule
$$H_\mu\triangleright H_\nu = L_{\mu\nu}^\si H_\si.$$
Equality (\ref{e1}) implies that the following
quasi-associativity property holds
\be
(H_\mu\triangleright H_\nu)\triangleright H_\si -
(H_\nu\triangleright H_\mu)\triangleright H_\si
=
H_\mu\triangleright (H_\nu\triangleright H_\si) -
H_\nu\triangleright (H_\mu\triangleright H_\si).
\label{e2}
\ee
Conversely, if a bilinear pairing $\triangleright$ on ${\bf H}$
satisfies this condition,
the skew-symmetric  operation
\be
[H_\mu,H_\nu]&=&(H_\mu\triangleright H_\nu) - (H_\nu\triangleright H_\mu)
\label{e3}
\ee
fulfils the Jackobi identity, and $\triangleright$ becomes
a left representation of the Lie algebra ${\bf H}$ equipped with
Lie bracket (\ref{e3}) on itself.

Lie algebra action $\triangleright$ induces  an action of the
Lie group ${\bf G}$ turning  ${\bf H}$
into the left ${\bf G}$-module. Consider now a 1-cocycle  $\varphi$ on
the group ${\bf G}$ with values in ${\bf H}$  \cite{G}.
This means  that $\varphi$ obeys the equation
\be
\varphi(xy)=\varphi(y)+ y^{-1}\triangleright\varphi(x),
\quad x,y \in {\bf G}.
\label{GrCoc}
\ee
Lie algebra  1-cocycle  $\partial\varphi$ is in one-to-one
correspondence with  $\varphi$, being its derivative taken
at the group identity \cite{G}. It satisfies the equation (cf. with
(\ref{e3}))
$
\partial\varphi([H_\mu,H_\nu])=H_\mu\triangleright\partial\varphi(H_\nu)-
                               H_\nu\triangleright\partial\varphi(H_\mu).
$
Suppose the linear mapping $\partial\varphi$ to be nondegenerate.
Then the identity map $id\colon{\bf H}\to {\bf H}$ is a 1-cocycle
with respect to the new action defined as
$(\partial\varphi)^{-1}\circ\triangleright \circ\partial\varphi$.
Thus, nondegenerate 1-cocycles of Lie algebras are in one-to-one
correspondence
with bilinear quasi-associative, in the sense of (\ref{e2}),
operations on ${\bf H}$.
Only non-degenerate cocycles are suitable for our purposes, so we
will think of them as of identity maps, and all the  freedom
will be encoded in the  choice of action $\triangleright$.
Note that a 1-coboundary normalized to $id$ implies the
existence of the right unity $H_e$ that is
$$
H_\mu\triangleright H_e = H_\mu.
$$
The group cocycle in terms of Lie algebra coordinates $\x^\mu$
in a neighbourhood  of the identity reads
$$\varphi^\mu(\x)=\Bigl(\frac{e^{-L(\x)}-1}{\scriptstyle-L(\x)}
\Bigr)^\mu_\nu\x^\nu$$
and the coboundary can be written as
$$\varphi^\mu(\x)=(1- e^{-L(\x)})^\mu_\nu\x_e^\nu,
\quad
H_\nu\x_e^\nu =H_e.
$$

Consider the semidirect sum ${\bf L}={\bf H}\triangleright {\bf H}^*$
with the Lie bracket given by (\ref{HHrel}) and (\ref{HXrel})
such that the condition (\ref{e1}) holds.
Since $\partial\varphi=id$ is non-degenerate,
the function $\varphi$ is invertible in a neighbourhood of the identity in
${\bf H}$.  Its inverse $\psi$ as well as $\varphi$ itself  are treated as
columns whose components are formal series in coordinate functions
generating ${\cal U}({\bf H}^*)$.
\newtheorem{theor}{Theorem}
\begin{theor}
   The element ${\cal F}=\exp(H_\nu \otimes  \psi^\nu(X))$ satisfies
the twist equation.
\end{theor}
\underline{Proof.} 
The element $\exp(H_\nu \otimes  \psi^\nu(X))$ satisfies the identity
(\ref{f-twist1}). If we prove the second identity (\ref{f-twist2}),
the theorem will be stated. 
Denote $\tX^\mu=\psi^\mu(X)$ and evaluate $\Delta_{\cal F}(X)$:
\be
 \Delta_{\cal F}(X^\mu)
 &=&\exp(H \otimes  \tX)( X^\mu\otimes  1+1\otimes  X^\mu )
 \exp(-H \otimes  \tX ) \n\\
 &=&  X^\nu\otimes (e^{-L(\tX)})^\mu_\nu+1\otimes  X^\mu .
\label{eq}
\ee
The map  $\Delta_{\cal F}(h) ={\cal F}\Delta(h){\cal F}^{-1}$
is an algebra homomorphism ${\cal U}({\bf L})\to 
{\cal U}({\bf L})^{\otimes 2}$.
Henceforth, (\ref{eq}) entails the equation
\be
\varphi(\Delta_{\cal F}(\tX))&=&e^{-L(1\otimes \tX)}\varphi(\tX\otimes  1)+
\varphi(1\otimes \tX) \,. 
\label{gr_coc}
\ee
Since $\varphi$ is invertible as a map of ${\bf H}$ on  ${\bf H}$, we
find $\Delta_{\cal F}(\tX^\mu) = {\cal D}^\mu(\tX\otimes 1, 1\otimes \tX)$
where ${\cal D}^\mu(\x_1,\x_2)$  is the Campbell-Hausdorf series.
This yields (\ref{f-twist2}) and therefore the twist equation
(\ref{TE}) for $\exp(H_\nu \otimes \psi^\nu(X))$ is valid. $\bullet$ 

Now we can evaluate the twisted
coproducts in  terms of new generators $\tX^\mu$.
A straightforward calculations show that
\begin{eqnarray}
\Delta_{\cal F}(H_\mu)    &=& H_\nu \otimes  g(\tX)^\nu_\mu +1\otimes  H_\mu,
    \label{HCopr}
\end{eqnarray}
where $g(\x)$ is a map ${\bf H}\to {\bf H}$ to be found.
Imposing coassociativity conditions we find that function $g$ realizes a
left group action on ${\bf H}$ which is generated by a Lie algebra
representation. To evaluate this action let us perform the following
Lie algebra isomorphism $H_\mu\to H_\mu$, $X^\mu\to \xi X^\mu$. The
specific form of the classical $r$-matrix allows us to consider $\xi$
as the deformation parameter.
Taking into account  $\frac{d}{d\xi}\tX^\mu(0)=X^\mu$, $\tX^\mu(0)=0$
and calculating  $\frac{d}{d\xi}{\cal F}\Delta(H_\nu){\cal F}^{-1}|_{\xi=0}$
we find
\be
\frac{d}{d\xi}\Delta_{\cal F}(H_\nu)|_{\xi=0}&=&
[H_\si \otimes  X^\si, H_\nu\otimes  1+ 1\otimes  H_\nu] \n\\
&=&
C^\mu_{\si\nu}H_\mu\otimes  X^\si + L^\mu_{\nu\si}H_\mu\otimes  X^\si=
L^\mu_{\si\nu}H_\mu\otimes  X^\si. \n
\ee
Performing this for the coproduct (\ref{HCopr}) and comparing the results
we find $g(\tX)=e^{L(\tX)}$.
Thus the coproduct on generators $H_\mu$, $\tX_\nu$ reads
\be
\Delta_{\cal F}(\tX^\mu)= {\cal D}^\mu(\tX\otimes  1,1\otimes  \tX) ,\quad
\Delta_{\cal F}(H_\mu)= H_\nu\otimes (e^{L(\tX)})^\nu_\mu  +1\otimes  H_\mu ,
\label{relHX0}
\ee
Using these relations it is easy to find also the antipodes,
\be
S_{\cal F}(\tX^\mu)=-\tX^\mu, \quad 
S_{\cal F}(H_\mu)=- H_\nu(e^{-L(\tX)})^\nu_\mu \,.
\label{inverse}
\ee
Expressing $\tX^\mu$ through $X^\nu$ we can evaluate
the twisted antipode on the classical generators as well.

\section{Conclusion}

The triangular deformation of the universal enveloping algebra 
of $sl(N)$ started already by Gerstenhaber et al \cite{GER}, 
was realized in this paper as a twisting with explicit form of the 
twist element $\cal F$ (\ref{twist-sl(N)}) (extended jordanian twist). 
The Hopf subalgebra of the type ${\cal U}_{\cal F}({\bf B}^{\vee})$ 
generated by the twist is self-dual. The twisted 
coproduct of the $sl(N)$ generators can be expressed as finite 
sums of classical generators and a function $\sigma$ of 
the highest root vector primitive with respect to the 
twisted coproduct. The commutation relations of the 
quantum space generators were defined using the twist $\cal F$ 
action on commutative coordinates. 
The cohomological properties of the involved Lie bialgebra 
permit to explain the connection of the Drinfeld-Jimbo (standard) 
quantization with this twisting. 

The explicit expression of the twist $\cal F$ gives rise to a possibility 
to evaluate the Clebsch-Gordan coefficients of the twisted $sl(N)$ 
in terms of the original CGC and the entries of the matrix 
$F = (\rho_1 \otimes \rho_2) {\cal F}$ \cite{KUL1}, as well as to get 
the relations among the FTR-approach generators $L^{(\pm)}$ 
of the twisted algebra and the generators of the original 
algebras. It can be 
used also to construct the quantum double \cite{DRIN,VLA}. 

The construction of the extended jordanian twist 
was generalized to certain class of 
inhomogenious Lie algebras, using properties of 
the Campbell-Hausdorf series. 
Further generalizations, in particular to 
Lie superalgebras \cite{KUL3}, twisting of the corresponding 
Yangians and new integrable models, twist elements for other  
boundary solutions to the classical Yang-Baxter 
equation \cite{GER2} are under study.

\end{document}